\documentclass[12pt]{amsart}

\usepackage{psfrag,epsfig}

\theoremstyle{plain}
\newtheorem{theorem}{Theorem}
\newtheorem{lemma}[theorem]{Lemma}
\newtheorem{corollary}[theorem]{Corollary}
\newtheorem{proposition}[theorem]{Proposition}
\theoremstyle{definition}

\newtheorem{remark}{Remark}

\def\Const{{\rm Const\ }}

\def\Lie{\mathop{\rm Lie}}

\def\EXP{{\mathbb{E}}}
\def\PROB{{\mathbb{P}}}

\def\integers{{\mathbb{Z}}}

\def\Int{{\rm Int}}

\def\diam{{\rm diam}}

\def\supp{{\rm supp}}
\def\l{\left}
\def\r{\right}

\def\bdef{\begin{Def}}
\def\endef{\end{Def}}
\def\bthm{\begin{theorem}}
\def\ethm{\end{theorem}}
\def\blm{\begin{lemma}}
\def\elm{\end{lemma}}
\def\bprop{\begin{proposition}}
\def\eprop{\end{proposition}}
\def\brm{\begin{remark}}
\def\erm{\end{remark}}
\def\beq{\begin{eqnarray}}
\def\eneq{\end{eqnarray}}
\def\beal{\begin{aligned}}
\def\enal{\end{aligned}}

\def\brx{{\bar x}}
\def\bry{{\bar y}}

\def\brlambda{{\bar\lambda}}

\def\bd{{\mathbf{d}}}

\def\cB{{\mathcal{B}}}
\def\cC{{\mathcal{C}}}
\def\cF{{\mathcal{F}}}

\def\cH{{\mathcal{H}}}

\def\cP{{\mathcal{P}}}
\def\cQ{{\mathcal{Q}}}
\def\cR{{\mathcal{R}}}
\def\cS{{\mathcal{S}}}
\def\cT{{\mathcal{T}}}
\def\cW{{\mathcal{W}}}

\def\R{\mathbb R}
\def\Z{\mathbb Z}
\def\T{\mathbb T}

\def\gm{\gamma}
\def\th{\theta}
\def\~{\tilde}

\def\lb{\lambda}

\def\vn{{\vec n}}

\author{D. Dolgopyat, V. Kaloshin and L. Koralov}
\thanks{D.D. was partially supported by NSF and Sloan Foundation,
V. K. was partially supported by American Institute of Mathematics
Fellowship and Courant Institute, and L. K. was partially supported
by NSF postdoctoral fellowship.}
\title{A Limit shape theorem for periodic stochastic dispersion.}

\begin{document}
\begin{abstract}
We consider the evolution of a connected set on the plane carried by
a periodic incompressible stochastic flow. While for almost every
realization of the random flow at time $t$ most of the particles are
at a distance of order $\sqrt{t}$ away from the origin, there is
a measure zero set of points, which escape to infinity at the linear
rate. We study the set of points visited by the original set by
time $t$, and show that such a set, when scaled down by the factor
of $t$, has a limiting non random shape.
\end{abstract}
\maketitle

\section{introduction.}

This paper deals with the long time behavior of a passive scalar
carried by an incompressible random flow. As has been demonstrated
for a large class of stochastic flows with zero mean, under some
mixing conditions on the flow, the displacement of a single
particle is typically of order $\sqrt{t}$ for large $t$. In
\cite{DKK} the authors show that for  almost every realization of
the random flow, if one considers the image of an open set under
the action of the flow, then its spatial  distribution, scaled by
the square root of time, converges weakly to a Gaussian
distribution. On the other hand, it has been shown in the work of
Cranston, Scheutzow, Steinsaltz, and Lisei ( \cite{CSS1},
\cite{CSS2}, \cite{LS}, \cite{SS}) that in any open  set there are
points which escape to infinity at the linear rate. Denote the
original set by $\Omega$. One can think of $\Omega$ as of an oil
spill or of a pollutant, say on the surface of the ocean. The
evolution of the set under the action of the flow will be denoted
by $\Omega_t$.

We shall study  the set of "poisoned" points, that is those
visited by the image of $\Omega$ before time $t$
\[
 \cW_t(\Omega)  = \bigcup_{s \leq t} \Omega_s~.
 \]
As shown in  \cite{CSS1} and \cite{CSS2} the diameter of this set
grows linearly in time almost surely. We shall be interested in
its limit shape.

Consider a stochastic flow of diffeomorphisms on $\R^2$ generated
by a finite-dimensional Brownian motion
\begin{equation}
\label{SF} dx_t=\sum_{k=1}^d X_k(x_t) \circ d\theta_k(t)+ X_0(x_t)
dt
\end{equation}
where $X_0, X_1, \dots , X_d$ are $C^\infty$ divergence free
periodic vector fields and $\vec\theta(t)=(\th_1(t), \dots,
\th_d(t))$ is a standard $\R^d$-valued Brownian motion.

We impose several  assumptions on the vector fields $X_0, X_1,
\dots , X_d$, which are stated in the next section (cf
\cite{DKK}). All those, except the assumption of zero drift,  are
{\it nondegeneracy assumptions} and are satisfied for a generic
set of vector fields $X_0, X_1, \dots , X_d$.

The main result of this paper is the following:

\begin{theorem}{\rm (The Shape Theorem)} \label{t1} Let the original
set $\Omega$ be bounded and contain a continuous curve. Under
Assumptions (A)-(E) from Section \ref{asm} on the vector fields,
there is a compact convex non random set $\cB$, independent of
$\Omega$, such that for any $\varepsilon> 0$ almost surely \beq
\label{shape} (1-\varepsilon)t \cB \subset \cW_t(\Omega) \subset
(1+\varepsilon)t \cB \eneq for all sufficiently large $t$.
\end{theorem}

In Section \ref{lb} we employ the estimates on the behavior of the
characteristic function of a measure carried by the flow in order
to bound from above the time it takes for the image of a curve to
reach a fixed neighborhood of a far away point (cf \cite{DKK}).
This bound in turn implies the lower bound in (\ref{shape}) for
the set $\cW_t$.

The key element in the proof of the upper bound for $\cW_t $ is to
show that the set $\cW_t$ for large $t$ is almost independent of
the original set (which, as will be demonstrated, can be taken to
be a curve). In order to prove this we show that given two bounded
curves $\gamma$ and $\gamma'$ we can almost surely find a contour,
which contains $\gamma'$ inside, and which consists of a finite
number of pieces of $\gamma_t$, and a finite number of stable
manifolds (whose length tends to zero as they evolve with the
flow). In this way we see that if a point is visited by the image
of $\gamma'$, then its small neighborhood is earlier visited by
the image of $\gamma$.

In Section \ref{cnt} we describe the construction of the contour,
and in Section \ref{ub} we complete the proof of the upper bound.

\section{nondegeneracy assumptions.}
\label{asm}  In this section we formulate a set of assumptions on
the vector fields, which in particular imply the Central Limit
Theorem for measures  and the estimates on the behavior of the
characteristic function of a measure carried by the flow (see
\cite{DKK}). Such estimates are essential in the proof of the
Shape Theorem. Recall that $X_0, X_1, \dots , X_d$ are assumed to
be periodic and divergence free.

(A) ({\it hypoellipticity for $x_t$}) For all $x\in \R^2$ we have
\[
{\Lie} (X_1,\dots,X_d)(x)=\R^2.
\]

Denote the diagonal in $\T^2 \times \T^2$ by
\[
\Delta=\{(x^1,x^2)\in \R^2 \times \R^2:\ x^1=x^2\
({\rm mod}~ 1) \}.
\]

(B) ({\it hypoellipticity for the two--point motion}) The
generator of the two--point motion $\{(x^1_t,x^2_t):\ t>0\}$ is
nondegenerate away from the diagonal $\Delta$, meaning that the
Lie brackets made out of $(X_1(x^1),X_1(x^2)),$ $\dots,
(X_d(x^1),X_d(x^2))$ generate $\R^2 \times \R^2$.

To formulate the next assumption we need additional notations. Let
$Dx_t:T_{x_0}\R^2 \to T_{x_t}\R^2$ be the linearization of $x_t$
at $t$. We need the hypoellipticity  of the process
$\{(x_t,Dx_t):\ t>0\}$. Denote by $TX_k$ the derivative of the
vector field $X_k$ thought as the map on $T\R^2$ and by
$S\R^2=\{v\in T\R^2:\ |v|=1\}$ the unit tangent bundle on $\R^2$.
If we denote by $\~X_k(v)$ the projection of $TX_k(v)$ onto
$T_vS\R^2$, then the stochastic flow (\ref{SF}) on $\R^2$ induces
a stochastic flow on the unit tangent bundle $S\R^2$ defined by
the following equation: \[ d\~x_t=\sum_{k=1}^d \~ X_k(\~x_t)\circ
d\th_k(t)+ \~X_0(\~x_t)dt. \] With these notations we have
condition

(C) ({\it hypoellipticity for $(x_t,Dx_t)$}) For all $v\in S\R^2$
we have
\[
\Lie(\~X_1,\dots,\~X_d)(v)=T_vS\R^2~.
\]

For measure-preserving stochastic flows with conditions (C)
Lyapunov exponents $\lb_1, \lb_2$  exist by {\it multiplicative
ergodic theorem for stochastic flows} of diffeomorphisms (see
\cite{C2}, thm. 2.1). Moreover, the sum of Lyapunov exponents
$\lb_1+ \lb_2$ should be zero (see e.g. \cite{BS}). Under
conditions (A)--(C) the leading Lyapunov exponent is positive
\[
\lb_1=\lim_{t\to \infty}
\frac 1t \log |Dx_t(x)(v)| > 0~.
\]
 Indeed, Theorem 6.8 of \cite{Bx} states that under
condition (A) the maximal Lyapunov exponent
$\lb_1$ can be zero only if for almost every
realization of the flow (\ref{SF}) one of the following two
conditions is satisfied

(a) there is a Riemannian metric $\bd$ invariant with respect to
the flow (\ref{SF}) or

(b) there is a direction field $v(x)$ on $\T^2$ invariant with
respect to the flow (\ref{SF}).

However (a) contradicts condition (B). Indeed, (a) implies that
all the Lie brackets of $\{(X_k(x^1),\ X_k(x^2))\}_k$
are tangent to the leaves of the foliation
\[
\{(x^1,x^2)\in \T^2 \times \T^2:\ \bd(x^1,x^2)=\Const\}
\]
and don't form the whole tangent space. On the other hand (b)
contradicts condition (C), since (b) implies that all the Lie
brackets  are tangent to the graph of $v$. This positivity of
$\lb_1$ is crucial for our approach.

Let $L_{X_k} X_k(x)$ denote the derivative of $X_k$ along $X_k$ at
the point $x$. Notice that $\sum_{k=1}^d L_{X_k} X_k + X_0$
is the deterministic components of the stochastic flow (\ref{SF})
rewritten in Ito's form. We require that the flow has no
deterministic drift, which
is expressed by the following conditions.

(D) ({\it zero drift})
\begin{equation}
\label{D1}
\int_{\T^2} \l(\sum_{k=1}^d L_{X_k} X_k + X_0\r)(x) dx =0~,
\end{equation}
\beq \label{D2}
\int_{\T^2} X_k(x) dx =0~,~~k=1,...,d~.
\eneq

Assumption (\ref{D2}) is not needed for the Shape Theorem to be true.
Still, in order to simplify the proof, that is in order to use
the results of \cite{DKK} without technical modifications, we shall
assume (\ref{D2}) to hold.

The last assumption is concerned with the geometry of the stream
lines of the vector fields $X_1,...,X_d$. Fix a coordinate system
on the $2$-torus $\T^2=\{x=(x_1,x_2)\ \textup{mod}\ 1\}$. As
the vector fields have zero mean and are divergence free, there are
periodic {\it stream} functions $H_1,...,H_d$, such that
$X_k(x) = (-H'_{x_2}, H'_{x_1})$. Also, by conition (A),
none of the
points $x \in \T^2$ is a critical point for all of the stream
functions $H_1,...,H_d$ at the same time.
We require the following

(E) ({\it condition on the critical points of $H_k$})
For each $k$ all of the
critical points of $H_k$ are non degenerate, and all the critical
values are distinct.

\section{Lower bound.}
\label{lb}
\subsection{An auxiliary statement.}
The following estimate will be repeatedly used in the proof of the
lower bound.
\begin{lemma}
\label{SMartLD} Let $\{\xi_j\}$ be a sequence of random variables
such that $$\EXP(\xi_{j+1}|\xi_1\dots \xi_j)\leq 0,$$ and that for
any $N$ the sequence $\{ \EXP|\xi_j|^N \}_j $ is bounded. Then for
any $\varepsilon > 0$ there exists $\kappa = \kappa(\varepsilon,
N)>0$ such that for each $n\in \Z_+$ we have
$$
\PROB \left\{ \sum_{j=1}^n \xi_j \geq  \varepsilon n \right\}
\leq \kappa n^{-N}.
$$
\end{lemma}
\proof  Define a sequence of random variables $\zeta_j =
\xi_j - \EXP(\xi_{j}|\xi_1\dots \xi_{j-1})$. Then $M_n =
\sum_{j=1}^n \zeta_j$ is a martingale whose  quadratic variation
is equal to $\langle M \rangle_n = \sum_{j=1}^n \zeta_j^2$. By the
martingale inequality
\[
\EXP M_n^{2N} \leq C_N \EXP \langle M \rangle_n^N \leq C'_N n^N~.
\]
Therefore, by the Chebyshev inequality
\[
\PROB \left\{ \sum_{j=1}^n \zeta_j \geq  \varepsilon n \right\}
\leq \PROB \left\{ M_n^{2N} \geq  (\varepsilon n)^{2N} \right\}
\leq  \kappa n^{-N}.
\]
Since $\xi_j \leq \zeta_j$ the result also holds
for the original sequence $\{\xi_j\}$. \qed

\subsection{Linear growth, an estimate from below.}
\label{eb} Let the initial set be a curve $\gamma_0$, and consider
a point $A$ in the plane. We shall estimate the tail of the
probability distribution of the time it takes for the curve to
reach  an $R$-neighborhood of $A$ in terms of the distance between
$\gamma_0$ and $A$ (the constant $R$ will be selected later). See
Figure 1.

We shall call a curve  {\it long} if its diameter is greater or
equal than 1. Without loss of generality we may assume that the
original curve is long. Given a curve $\gamma_0$ and a point  $A$
we define $\tau^R = \tau^R(\gamma_0, A)$ to be the first moment of
time when the image of $\gamma_0$ reaches the $R$-neighborhood of
$A$, and at the same time the image of $\gamma_0$ is long, that is
\beq \label{stop}
\tau^R = \inf \{t>0 :\ {\rm dist} (\gamma_t,A) \leq R, ~
\diam(\gamma_t) \geq 1 \}~.
\eneq

Let  $\cW_t^R(\gm_0)$ be the $R$-neighborhood of $\cW_t(\gm_0)$,
that is the set of points, whose $R$-neighborhood is visited by
the image of the original set before time $t$,
\[
\cW_t^R(\gm_0) = \{x: {\rm dist}(x, \gamma_s) \leq R~~~~{\rm
for}~~ {\rm some}~~s \leq t \}~.
\]

\bprop \label{LmP-PLD} Let $d = {\rm dist}(A, \gamma_0)$.
There is  a constant $R > 0$, such that for any
$m$ there exist positive $C$ and $ \beta$ such that for any long
curve $\gamma_0 \subset \R^2$ and any point $A\in \R^2$ we have
$$
\PROB\left\{\tau^R>\beta d \right\}\leq C d^{-m}.
$$
\eprop

The proof of Proposition \ref{LmP-PLD}  will rely on two other
lemmas. We first state the lemma (proved in \cite{CSS2}), which
shows that $\gamma_t$ can not grow faster than linearly.
\begin{lemma}
\label{NotTooFast} (\cite{CSS2}) Let  $\gamma_0$
be the initial curve and
$$
\varPhi_t=\sup_{x_0\in\gamma_0}
\sup_{0\leq s\leq t} ||x_s-x_0||.
$$

(a) There is a constant $C$ such that almost surely
$$
\lim\sup_{t\to\infty} \frac{\varPhi_t}{t}\leq C,
$$

(b) For any positive $r$ and $\alpha$ we have
$$ \sup_{t \geq 1}
\EXP\left(\exp\left[ \frac{r \varPhi_t^2}{t^2 \max(1,\ln
(\varPhi_t/t))^{2+\alpha}} \right]\right) <\infty.
$$
\end{lemma}

\begin{figure}[htbp]
\label{Cone}
  \begin{center}
    \begin{psfrags}
     \psfrag{gm}{$\gm$}
     \includegraphics[height=3.5in, width= 4in,angle=0]{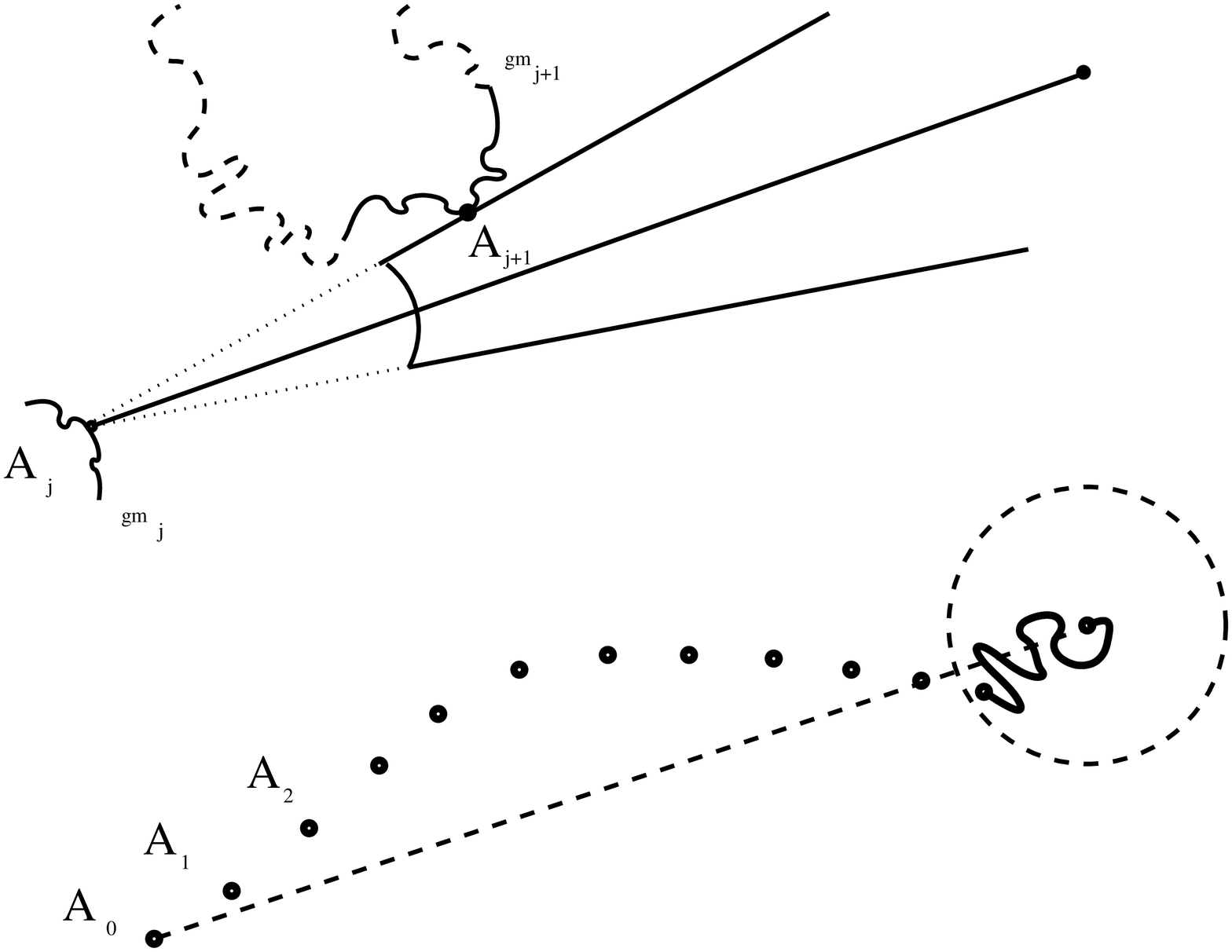}
    \end{psfrags}
    \caption{}
  \end{center}
\end{figure}

Given a pair $(A_0,\gamma_0)$, where $\gamma_0$ is a long curve,
and $A_0$ is a point on $\gamma_0$ we define inductively the
sequence $(A_j,\gamma_j)$ as follows. Suppose that $(A_{j},
\gamma_{j})$ is defined so that
\begin{itemize}
\item $\gamma_j$ is a connected interval of the image
of $\gamma_0,$ i.e. $\gamma_j\subset\gamma_{t_j};$

\item $A_j\in\gamma_j$;

\item $\gamma_j$ is long.
\end{itemize}
Given $\alpha\in(0,\frac{\pi}{2})$ define the truncated
$\alpha$-cone
$$
K_j(\alpha)=\{x: {\rm
dist}(x,A_j)\geq 1\quad {\rm and}\quad \angle (xA_jA) \leq
\alpha\},
$$
where $\angle (xA_jA)$ is the angle between the
segments $[A_j, x]$ and $[A_j, A]$.

Let $\gamma_j(t)$ be the image of $\gamma_j$ under the flow (so
that $\gamma_j=\gamma_j(t_j).$) Let $t_{j+1}$ be the first moment
such that
\begin{itemize}
\item $t_{j+1} - t_j \geq 1$;

\item $\diam(\gamma_j(t_{j+1}))\geq 1$;

\item $\gamma_j(t_{j+1})\bigcap K_j(\alpha)\neq\emptyset.$
\end{itemize}

Let $A_{j+1}$ be an arbitrary point in $\gamma_j(t_{j+1})\bigcap
K_j(\alpha)$, $B_R(A_{j+1})$ denote the $R$-ball around $A_{j+1}$,
and $\gamma_{j+1}$ be a connected interval of
$\gamma_j(t_{j+1})\bigcap B_1(A_{j+1})$ containing $A_{j+1}.$ Let
$\cF_{t_j}$ be the $\sigma$ algebra of events determined prior to
time $t_j$.

\begin{lemma}
\label{SectWait} There exists  $ \alpha<\frac{\pi}{2}$ such that
for any $ m$  we have
\[
\EXP((t_{j+1}-t_j)^m)| \cF_{t_j})<C_m.
\]
\end{lemma}

\proof It is sufficient to prove that $\EXP (t_1^m) < C_m$ with
$C_m$ independent of the original long curve $\gamma_0$. Without
loss of generality we may assume that the point $A_0$ on the curve
$\gamma_0$ coincides with the origin, and that the  original curve
is contained in a ball of radius 2 centered around the origin.

Fix some $0 <p < 1$. Since $\gamma_0$ is long, the Frostman Lemma
(\cite{Ma}) implies that there exists a probability measure $\nu$
on $\gamma_0$ whose $p$ energy is bounded
\begin{equation}
\label{st1}
 I_p(\nu) = \int\int_{\gamma_0\times\gamma_0} \frac{
d\nu(x) d \nu(y)}{|x-y|^p} \leq C_p ~.
\end{equation}

Let $f(x)$ be a continuous non-negative function compactly
supported inside $K_0(\alpha)$. By the results of \cite{DKK} under
condition (\ref{st1})  there exists a non degenerate $2\times 2$
matrix $D$ such that for any $\rho, m>0$ there exists $T$ such
that for all $t>T$

\begin{equation}
\label{MesSect} \PROB \left\{ |\int_{\gamma_0}  f \left(
\frac{x_t}{\sqrt{t}} \right) d \nu - \bar f |> \rho \right\}
\leq t^{-m}~,
\end{equation}
where  $\bar f$ denotes the integral of $f$ with respect to the
Gaussian measure $\nu$ with zero mean and variance $D.$ The last
inequality of Sec. 7 in \cite{DKK} establishes (\ref{MesSect}) for
functions of the form $f(x)=\exp(i\xi x)$, but this inequality is
also valid for compactly supported functions, since they could be
uniformly approximated linear combinations of complex
exponentials. Thus, for large $t$

$$ \PROB\left\{\gamma_t\bigcap \sqrt{t}\ \supp(f)\ = \emptyset
\right\}\leq t^{-m}. $$
Since for $t\geq 1$, we have that $\sqrt{t}\
\supp(f)  \subset K_0(\alpha)$ we get that for $t\geq T$
$$
\PROB\{\gamma_t\bigcap K_0(\alpha)= \emptyset \} \leq t^{-m}.
$$
From (\ref{MesSect}) it also follows that for large $t$
\[
\PROB \{ {\rm diam}(\gamma_t) < 1 \} \leq t^{-m}~.
\]
 Since $m$ is arbitrary, this implies the required
result. \qed

{\it Proof of Proposition \ref{LmP-PLD}.} Let $r_j={\rm dist}(A_j, A)$.
Let us show that there exists $ R>0$ such that for any $ m$ there
is $ \beta$ for which
\beq \label{f1}
\PROB \{r_j>R\quad{\rm
for}\quad j=1,2,\dots ,[\beta r_0]\} \leq \frac{C_m}{r_0^m}~.
\eneq
 We first establish a weaker bound: there is $R>0$ such that
\beq \label{wb}
\PROB \{r_j > R \quad {\rm for} \quad
 j=1,2,\dots ,[r_0^m])\leq \frac{C_m}{r_0^{m-1}}~.
\eneq
Note that there exist positive constants $R$ and $K$ such
that if $r_0 > R$, then $\EXP (r_1 - r_0) \leq -K$. Let $j_1$ be
the first moment when $r_j \leq R$, and in general let $\{j_n\}_n$
be the times of the consecutive visits of $A_j$ to $B_R(A)$,
that is $A_{j_n} \in B_R(A)$.

 Then
\beq \label{sm}
r_{\min(j,j_1)}+K\min(j,j_1)
\eneq
 is a supermartingale. Hence
$$
r_0\geq \EXP(r_{\min([r_0^m],j_1)}+K\min([r_0^m],j_1))
\geq K[r_0^m]\ \PROB\{j_1>[r_0^m]\}~,
$$
which implies (\ref{wb}).

Recall the definition of the points $A_j$. Notice that for
some constants $C_1, C_2$ and $R_1$, which depend on $\alpha$,
we have if dist$(A_{j+1},A_j)\leq C_2 r_j$, then $r_{j+1}<r_j-R_1$.
Indeed, $A_{j+1}$ belongs to the truncated cone $K_j(\alpha)$
with the origin at $A_j$ and directed toward $A$. This implies
that the sequence $r_0, r_1, \dots , r_{[C_1r_0]}$  contains
an element $r_j$ such that either $r_j \leq r_0^{\frac{1}{m}}$ or
${\rm dist}(A_{j+1},A_j) > C_2 r_0^{\frac{1}{m}}$.

Let now $j^*$ be the first moment when $r_j \leq
r_0^{\frac{1}{m}}.$ We have
$$
\PROB\{j_1 \geq 2 C_1 r_0\}\leq
\PROB\{j^*> C_1 r_0\}+
\PROB\{j_1\geq 2 C_1 r_0|j^*\leq C_1 r_0\}
= P_1 + P_2.
$$

In order to estimate the first term on the right hand side of this
inequality we note that if $j^*>C_1 r_0$ then ${\rm
dist}(A_{j+1},A_j)> C_2 r_0^{\frac{1}{m}}$ for some $j \leq [C_1
r_0].$ From  Lemmas  \ref{NotTooFast} and \ref{SectWait} it now
follows that $P_1$ decays faster than any power of $r_0$, in
particular $P_1 \leq \frac{C_3}{\sqrt{r_0}}$.

We estimate the second term on the right hand side using (\ref{wb})
with $m = 2$
\[
P_2 \leq C_3 r_0^{-\frac{m-1}{m}} = \frac{C_3}{\sqrt{r_0}}~.
\]

Therefore,
\beq \label{kt} \PROB\{j_1 \geq 2 C_1
r_0\}\leq\frac{C_4}{\sqrt{r_0}}~.
\eneq
 Let a positive integer $N$ be fixed. Then
\[
 \PROB\{j_1 \geq 2 C_1 N r_0\} =
\PROB\{j_1 \geq 2 C_1 N r_0 ;\ r_{[2 C_1 r_0]} > r_0\}  +
 \]
 \[
  \PROB \{j_1 \geq 2 C_1 r_0 ;\ r_{[2 C_1 r_0]} \leq r_0\}\
 \PROB \{j_1 \geq 2 C_1 N r_0 |\ j_1
\geq 2C_1 r_0 ; r_{[2 C_1 r_0]} \leq r_0\}
\]
\[ = Q_1 + Q_2 Q_3~.
\]
By Lemmas  \ref{NotTooFast} and \ref{SectWait}  the tail of the
probability distribution of $r_{j+1} -r_j$ decays faster than any
power. From this fact and  the fact that (\ref{sm}) is a
supermartingale by Lemma \ref{SMartLD} it is seen that $Q_1 $
decays faster than any power of $r_0$.

The factor $Q_2$ on the right hand side of the last inequality is
estimated from above by $\frac{C_4}{\sqrt{r_0}}$ due to (\ref{kt}).
Using the Markov property and continuing by induction we see that
$Q_3$ is estimated from above by $C_5 r_0^{-\frac{N-1}{2}}$.
We conclude that
\[
 \PROB\{j_1 \geq 2 C_1 N r_0\} \leq C_6 r_0^{-\frac{N}{2}}~.
\]
 Since $N$ was arbitrary this
implies (\ref{f1}).

Note that the estimate (\ref{f1}) remains valid even if the
distance from the original point $A_0$ to $A$ is less than $r_0$.
Therefore  any $m>0$ there are constants $C_m $ and $\beta$ such
that
\[
\PROB \{t_{j_1} > \beta d \} \leq \frac{C_m}{d^m}
\]
for all $d$. By the Markov property
\begin{equation}
\label{betweenvisits}
\EXP(t_{j_{n+1}} - t_{j_n} | \cF_{t_{j_n}} ) \leq C
\end{equation}
for some $C$ and for all $n$. We can therefore apply Lemma
\ref{SMartLD} with $\xi_n = t_{j_{n+1}} - t_{j_n} - C$ to obtain
that for any $N$ there exists $\kappa$ such that \beq \label{tme}
\PROB\{t_{j_n}>2Cn\}\leq  \kappa   n^{-N} . \eneq This combined
with (\ref{f1}) completes the proof of the proposition. \qed

\begin{corollary} \label{cor1} There exist a positive constants
$c $ and $R$ such that almost surely \\ $  B_{ct}(0) \subset
\cW_t^R$ for large $t.$
\end{corollary}

\proof Consider a covering of $B_{ct}(0)$ by balls of radius
$R/2$. By Proposition \ref{LmP-PLD} for each of the  balls of
radius $R/2$, the probability that it is not visited by the curve
by time $t$ decays exponentially in $t$, provided that $c$ is
small enough, and $R$ is large enough. On the other hand, for each
$c$ and $R$ the number of balls needed to cover $B_{ct}(0)$ grows
like $t^2$ times a constant. Therefore, the probability that the
$R$ neighborhood of some point in $B_{ct}(0)$ is not visited by
the curve before time $t$ decays exponentially in $t$. The
Corollary follows by Borel-Cantelli Lemma. \qed

From now on we fix  $R$  for which  Proposition \ref{LmP-PLD} and
Corollary \ref{cor1} hold.

Note that the bounds we obtained in the proof of
Proposition \ref{LmP-PLD} are uniform over all long curves. Let us
employ this fact in the following corollary. Let $\cC$ be the family
of long curves which lie completely inside $B_{2R}(0)$
(we may assume that $R >1$).
\begin{corollary}
\label{UI} The family of stopping times  $$ \left\{
\frac{\tau^R(\gamma_0,tv)}{t}\right\}_{t\geq 1, ||v||=1,
\gamma_0\in\cC} $$ is uniformly integrable.
\end{corollary}

\subsection{Stable norm.}
We shall now use the asymptotic of  $\tau^R$ (the time it takes a
curve to reach the $R$-neighborhood of a far away  point) in order
to define the limiting shape $\cB$.  Consider $$ |v|^R =\sup_{\cC}
\EXP\  \tau^R(\gamma,v). $$

Due to the stationarity of the underlying Brownian motion, and
due to the periodicity of the vector fields we have
 $$
\EXP\  \tau^{2R}(\gamma,(t_1+t_2) v)\leq \EXP\  \tau^R(\gamma, t_1
v)+ \EXP\  \tau^{R}(\gamma_1, t_2 v)~, $$ where $\gamma_1 \in \cC$
is some integer translation of a part of $\gamma_{\tau(\gamma,t_1
v)}$. Since by Proposition \ref{LmP-PLD}
\[
\EXP\  \tau^{R}(\gamma, (t_1+t_2) v) \leq \EXP\  \tau^{2R}(\gamma,
(t_1+t_2) v)+ C
\]
for some $C>0$, it follows that the function $|t v|^R + C $ is
sub-additive. Let
$\|v\|^R=\lim_{t\to\infty} \frac{|t v|^R}{ t} .$
 Similarly for $0\leq s\leq 1$
$$ |t( s v_1+(1-s) v_2) |^R \leq |t s v_1|^R +|t (1-s) v_2|^R+ C,
$$ so $$ \|s v_1+(1-s) v_2 \|^R \leq s \|v_1\|^R +(1-s)\|v_2\|^R .
$$ Let $\cB=\{v:\|v\|^R \leq 1\}.$ By the remarks above and by
Corollary \ref{UI} we have that $\cB $ is a convex compact set. It
will be seen that the norm $\|v\|^R$ and the set $\cB$ are
independent of $R$.

\begin{lemma}
\label{Contains} For any  $ \varepsilon > 0 $ almost surely
$(1-\varepsilon)t\cB \subset \cW_t^R $ for large enough $t$.
\end{lemma}

\proof It suffices to show that for any $ v$ and $  N$
$$
\PROB\{\tau^R(\gamma,tv)\geq (1+\varepsilon) t \|v\|^R\}\leq C_N
t^{-N}.
$$
By the definition of $\|v\|^R$ there exists $t_0$ such that

\beq \label{es1}  \EXP\ \tau^R(\gamma, tv)\leq
t(1+\frac{\varepsilon}{3}) \|v\|^R \eneq

for any $t\geq t_0$ and $\gamma \in \cC $. Define the stopping
time $\tau_1$ as follows
\[
\tau_1 = \inf\{ t: \gamma_t \bigcap B_R(t_0v) \neq \emptyset ;~
{\rm diam}(\gamma_t) \geq 1 \}~.
\]

Let $\gamma^{(1)}$ be a part of $\gamma_{\tau_1}$, which is long,
is contained in $B_{2R}(t_0 v)$,  and has a non empty intersection
with $B_R(t_0v)$. Similarly we define $\tau_2$ to be the first
time following $\tau_1$, when the image of $\gamma^{(1)}$ is long
and intersects $B_R(2 t_0 v)$, and let $\gamma^{(2)}$ be a part of
the image of $\gamma^{(1)}$, and so on. We have therefore
constructed a sequence of stopping times, such that
 $$
\tau^R(\gamma, n t_0 v)\leq \sum_{j=1}^n (\tau_j- \tau_{j-1})~.
$$
By (\ref{es1}), due to the periodicity of the underlying vector
fields, for large enough $t_0$ we have
\[
\EXP\ (\tau_j- \tau_{j-1})
\leq t_0\l(1+\frac{2 \varepsilon}{3}\r) \|v\|^R.
\]
Now the result follows by Lemma \ref{SMartLD} and Proposition
\ref{LmP-PLD}.\qed

Now we prove that Lemma \ref{Contains} remains valid even when
the $R$-ne\-i\-gh\-borhood of $\cW_t$ is replaced by $\cW_t(\gm)$ itself.

\begin{theorem}
\label{Contains2} For any  $ \varepsilon > 0 $ we have almost
surely $(1-\varepsilon)t\cB \subset \cW_t(\gm) $ for large enough
$t$.
\end{theorem}

This Theorem is a consequence of Proposition \ref{LmP-PLD} and the
fact that when a long curve reaches an $R$-neighborhood of a
point, the distribution of the time it then takes for the  curve
to sweep the entire neighborhood has fast decreasing tail. Thus
Theorem \ref{Contains2} follows from the standard Borel-Cantelli
arguments and the following

\begin{lemma} \label{lem}
 Let $\gamma_0$ be a long curve such that
${\rm dist}(\gamma_0, A)\leq R.$  Let
\[
\sigma = \inf\{t>0: B_R(A)
\subset \bigcup_{s \leq t} \gamma_s;
~ {\rm diam} (\gamma_t) \geq 1 \}~.
\]
Then for any $m>0 $ and some $C_m$ we have
\beq \label{es3}
\PROB(\sigma>t)\leq C_m t^{-m}~.
\eneq
\end{lemma}

\begin{figure}[htbp]
\label{Circle}
  \begin{center}
    \begin{psfrags}
     \psfrag{gh}{$\gamma_{1/2}$}
     \psfrag{g1}{$\gamma_1$}
     \psfrag{Bd}{$B_\delta(\tilde A)$}
     \includegraphics[height=3.5in, width= 4in,angle=0]{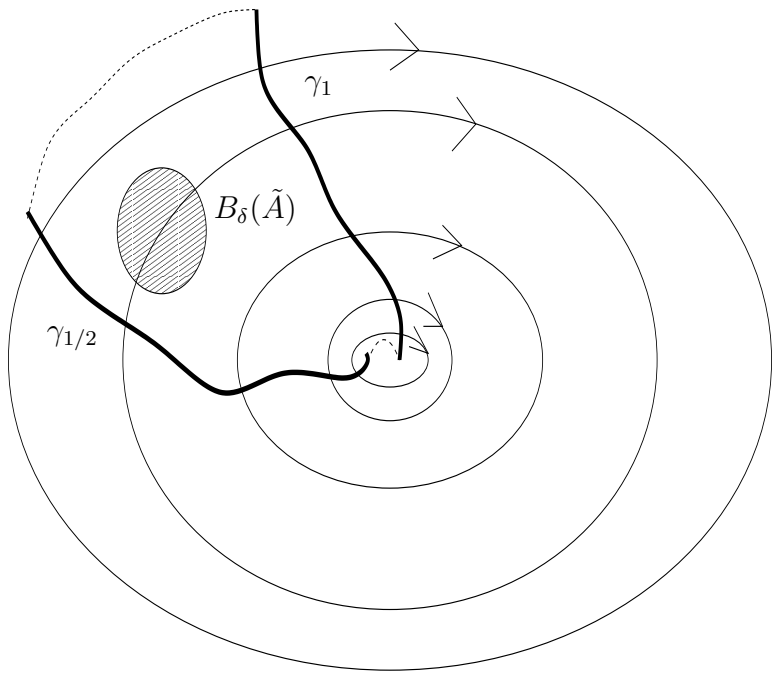}
    \end{psfrags}
    \caption{The image of a sufficiently long curve $\gamma$ covers
             a small ball $B_\delta(\tilde A)$ with positive probability.
             \newline
             (Dashed lines on the picture show the orbits of the
             endpoints of $\gamma$.)}
  \end{center}
\end{figure}

\proof  Since the ball $B_R(A)$ can be covered by a finite number
of arbitrarily small balls, it is sufficient to demonstrate
(\ref{es3}) for
\[
 \sigma_1 = \inf\{t: B_{\delta} (\tilde A) \subset
\bigcup_{s \leq t} \gamma_s;
~ {\rm diam} (\gamma_t) \geq 1 \}~,
\]
where $\tilde A \in B_R(A)$. As before,  let $j_n$ be the times of the
consecutive visits of $A_j$ to $B_R(A)$. We claim that for small
enough $\delta$  there exists $\theta<1$ such that \beq
\label{forml1}  \PROB\{\sigma_1 \geq t_{j_n}|\sigma_1 \geq
t_{j_{n-1}}\}\leq \theta~. \eneq Since $t_{j_n}-t_{j_{n-1}} \geq
1$ it is sufficient to show that any point $\tilde A$ with probability
$1-\theta$ has a neighborhood, which is swept by the time 1
evolution of any long curve $\gamma_0$, which initially is at a
distance no greater than $2R$ away from $\tilde A$.

Due to condition (A) we may assume that $\tilde A$ is not a critical
point of say $X_1$.  Let $\Gamma$ be the streamline of $X_1$, which
passes through $\tilde A$. Then, due to the periodicity, $\Gamma$ is
either a closed curve, or an unbounded periodic curve, or $\tilde A$
belongs to a separatrix (level set containing a saddle point) of
$H_1$.

 Consider the case when $\Gamma$ is a closed curve. Let
$U^{2\delta}$ be the region formed by the streamlines of $X_1$
passing through the $2\delta$ neighborhood of $\tilde A$. For
$\delta$ small enough the region $U^{2\delta}$ is homeomorphic to
an annulus. Let $B^1$ and $B^2$ be the interior and the exterior
circles bounding $U^{2\delta}$, and let $P$ be a point inside
$B^1$. Define the winding number of a point $x$ on a curve
$\gamma$ as the change of argument of rotation around $P$ for the
part of $\gamma$ from the initial point to $x$.  Due to the
hypoellipticity of the two point motion with probability $p_1
> 0$ we have the following: one of the points on
$\gamma_{\frac{1}{2}}$ lies inside $B^1$, while another point on
$\gamma_{\frac{1}{2}}$ lies outside $B^2$, and the winding number
around  $P$ is bounded in absolute value by a constant $K$ for all
points of $\gamma_{\frac{1}{2}}$. Notice, that in this case
$\gamma_{\frac{1}{2}}$ crosses each of the streamlines, which form
$U^{2\delta}$.

Let $\phi^1_t(x)$ be the evolution of $x$ under the action of the
vector field $X_1$: \[ \dot{\phi^1_t}(x) = X_1(\phi^1_t(x))~,~~~~
\phi_0(x) = x~. \] Let $T$ be the supremum over $U^{2\delta}$ of
the time it takes for a point to make a complete rotation
\[
T = \sup_{x \in U^{2\delta}} \inf_{t > 0} \{ t: \phi^1_t(x) = x \}~.
\]
Given the event that $\gamma_{\frac{1}{2}}$ crosses each of the
streamlines of $U^{2\delta}$,  for any $\varepsilon
>0$ there is a positive probability $p_2
>0$ that the time $t = \frac{1}{2}$ dynamics of (\ref{SF}) is
$\varepsilon$-close to the solution of
 \beq \label{forml2} \dot{\phi^1_t}(x) =
4T(K+1) X_1(\phi^1_t(x))~,~~~~ \phi^1_0(x) = x~. \eneq By
selecting $\varepsilon$ small enough we see that with positive
probability $p_2$ the time $\frac{1}{2}$ evolution of
$\gamma_{\frac{1}{2}}$ sweeps the entire annulus $U^{\delta}$.
Therefore, the probability that $B_{\delta}(\tilde A)$ is swept by
the image of $\gamma_0$ before $t = 1$ is not less than $p_1p_2 >
0$. We can take $\theta = 1- p_1p_2$. The cases when $\Gamma$ is a
periodic unbounded curve, and when $\tilde A$ belongs to a
separatrix of $H_1$ are treated similarly. (In the case when
$\tilde A$ belongs to the separatrix we may need to use one of
the vector fields $X_2,...,X_d$, or one of the commutators in
order to show that the evolution of $\gamma_0$ may pass through a
saddle point).

From (\ref{forml1}) it follows by induction that
\begin{equation}
\label{IndClose}
\PROB\{\sigma_1 \geq t_{j_n}\}\leq \theta^n
\end{equation}
for some $\theta<1.$   This combined with estimate (\ref{tme})
completes the proof of the lemma. \qed \\

\subsection{Occupation times.}
The next result is used in Subsection \ref{SSTSL}.

\begin{lemma} \label{occupation}
For any $\varepsilon>0$ there is a constant $R^*$ such
that for any curve $\gamma$ for any $A$
$$
\liminf_{T\to\infty} \frac{|  \{t\leq T:
\gamma_t\bigcap B_{R^*}(A)\neq\emptyset \} |}{T}
\geq 1-\varepsilon.
$$
\end{lemma}

\proof Let $j_n$ be as in Proposition \ref{LmP-PLD}. Let
$$
a_n(R^*)=\left|
\{t>0: t_{j_{n-1}} \leq t < t_{j_n} \ \ \ \textup{ and}\
\ \ A_{j_{n-1}}(t-t_{j_{n-1}})\not\in B_{R^*}(A)\} \right|
$$
be the amount of time between two consecutive visits of
a smaller $R$-ball $B_R(A)$ spend outside a bigger
$R^*$-ball $B_{R^*}(A)$. Then finiteness of expectation
of time (\ref{betweenvisits}) between consecutive visits
of $B_R(A)$ implies that for each $\varepsilon$ there is a
sufficiently large $R^*$ such that for all $n$
$$ \EXP\
a_n(R^*) <\frac{\varepsilon}{2}.
$$
Using Lemma \ref{SMartLD} we
obtain that almost surely $ \sum_{n=1}^N a_n(R^*)< \varepsilon N $
for a large $N.$ On the other hand by the definition of $t_j$
$$\sum_{n=1}^k \left(t_{j_n}-t_{j_{n-1}}\right)>k. $$ This implies
the required result. \qed

\begin{corollary}
\label{Simult}
There is a constant $R^*$ such that for any two open sets
$U_1, U_2$, any two points $A_1, A_2$, and some $T$ almost
surely there is a moment $t>T$ such that
$\phi_t U_j$ intersects $B_{R^*}(A_j)$ for both $j=1,2$.
\end{corollary}

\section{Geometry of stable manifolds.}
\label{cnt}
\subsection{Transitivity of stable lamination.}
\label{SSTSL} Before we can proceed further we need to establish
some properties of unstable lamination. Recall that one of our
assumptions is $\lambda_1>0.$ Let $0< \brlambda_1<\lambda_1$ then
by stable manifold theorem \cite{Ki} for all $t$ almost surely for
almost every $x$ the set $$ W^s(x,t)=\{y:\ \textup{for some}\
C(y)\ \textup{we have}\ d(f_{t,s} y, f_{t,s} x)\leq C(y)
e^{-\brlambda_1 (s-t)}\} $$ is a smooth curve passing through $x.$
By stationarity the parameters of $W^s$ such as its length,
curvature etc. have distributions independent of $t.$ Here we
prove the following.
\begin{theorem}
\label{WsTr}
Almost surely $W^s(\cdot,t)$ is transitive for all $t,$ that is,
given two open sets $U_1$ and $U_2$ for some $x\in U_1$
we have $W^s(x,t)\bigcap U_2\neq \emptyset$ or
$U_1$ and $U_2$ are connected by a stable manifold.
\end{theorem}
\proof Since $W^s(\cdot,t)=\varphi_t W^s(\cdot, 0)$ and the image
of a transitive lamination is transitive, it is enough to show
that $W^s(\cdot, 0)$ is transitive. Let $E_{U_1, U_2}$
denote the event
$$
\{\{\th(t)\}_{t>0}:\ \ \textup{for some}\ y\
W^s(y,0)\bigcap U_1\neq \emptyset \ \ \ \textup{and} \ \ \
W^s(y,0)\bigcap U_2\neq \emptyset\}
$$

We want to show that $\PROB\{E_{U_1, U_2}\}=1.$ In fact a weaker
statement suffices.

\begin{proposition}
\label{Alwposs}
Let $\cF_t$ denote the $\sigma$-algebra generated by
$\{\theta_j(s)| s\leq t\}$
Then for any $t>0$ there is a constant $\delta>0$ such that
$$
\PROB\{E_{U_1, U_2}|\cF_t\}\geq \delta .
$$
\end{proposition}

The rest of this subsection is devoted to the proof of Proposition
\ref{Alwposs}. This proposition implies the Theorem because
$\PROB\{E_{U_1, U_2}|\cF_t\}$ converges as $t\to+\infty$ to the
indicator of $E_{U_1,U_2}$ almost surely and so we get that this
indicator is greater than $\delta/2$ almost surely. Since the
indicator assumes only values 0 and 1 Theorem \ref{WsTr} follows.
\qed
\begin{lemma}
\label{CouldCr}
For any point $x\in \R^2$, positive $r_0,\ \varepsilon$ there exists
$p_1>0$ such that for all $y\in\R^2$ such that $d(x,y)\leq r_0$ we have
$$
\PROB\{W^s(x,0)\bigcap B(y,\varepsilon)\neq\emptyset\}>p_1.
$$
\end{lemma}

\proof There are positive $r_1, p_1$ such that for each pair:
a point $x\in \R^2$ and a time $t>0$ we have
$\PROB\{\diam(W^s(x,t)) > r_1\}>p_1$ or probability to have a stable
manifold at a point at a certain time of diameter at least a positive
$r_1$. Thus for each positive $\varepsilon_1$ there is a constant
$p_2>0$ such that for $x\in \R^2$ there is a point $y\in \R^2$ such
that $d(x,y)=r_1/2$ and for all $t$ we have
$\PROB\{W^s(x,t)\bigcap B(y,\varepsilon_1)\neq \emptyset\}>p_2.$
Another way to put it is that if a stable manifold at a point at
time $t$ is of diameter at least $r_1$, then it intersects a small
ball centered at a nearby point with a positive probability. Now by
hypoellipticity of the two-point motion there are positive
$\varepsilon_1$ and $p_3$ such that for each pair of points
$\brx, \bry \in \R^2$ and $d(x,y)\leq r_0$ for the time $1$ map
of the flow (\ref{SF}) we have
$$
\PROB\{\varphi_1 B(y,\varepsilon)\supset
B(\bry,\varepsilon_1)|x_1=\brx\}>p_3 .
$$
However, $W^s(x_1,1)$ crosses $\varphi_1 \left(B(y,\varepsilon)\right)$
with positive probability, therefore,
$$
\PROB
\l\{W^s(x,t)\bigcap B(y,\varepsilon)\neq\emptyset\r\}>p_2 p_3>0 .
$$
This completes the proof. \qed

\begin{figure}[htbp]
  \begin{center}
    \begin{psfrags}
      \psfrag{phi}{$\varphi_{t, t+\varepsilon} \Lambda $}
      \psfrag{q}{$q$}
      \psfrag{Q1}{$\cQ_1$}
      \psfrag{Q2}{$\cQ_2$}
        \includegraphics[height=2in, width= 4.5in,angle=0]{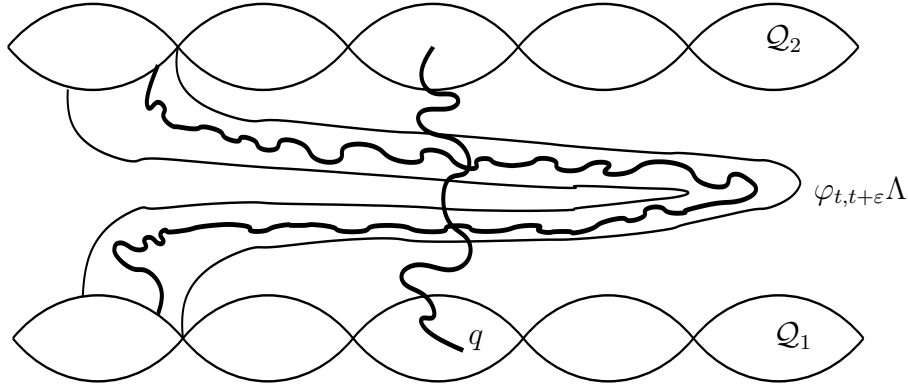}
    \end{psfrags}
    \caption{A stable manifold joining $\cQ_1$ and $\cQ_2$ should intersect
             the image of $U_1$ unless it is long enough to avoid the
             tongue $\varphi_{t, t+\epsilon} \Lambda$}
  \end{center}
\end{figure}

{\it Proof of Proposition \ref{Alwposs}.} Fix $1\leq k \leq d$ and
consider the vector field $X_k$ along with its stream function
$H_k:\T^2 \to \R$. By condition (E) $H_k$ exists and its lift
$\hat H_k:\R^2 \to \R$ has unbounded level sets $\hat
H_k^{-1}(c^*)$ for some regular values $c^*$. Choose two adjacent
critical values $c^-, c^+$ of  $H_k$ with $c^-<c^*<c^+$. Then the
stip of level sets $\cH_k=\{\hat H_k^{-1}(c^-,c^+)\}$ consists of
only unbounded level sets of $\hat H_k$ which are also stream
lines of $X_k$. Let $\cQ_1$ and $\cQ_2$ be the components
consisting of bounded streamlines bounding $\cH_k$ (see figure 3).
Choose $R^*>0$ from Corollary \ref{Simult}. Combining Lemma
\ref{occupation} and Corollary \ref{Simult} we have that given two
point $A_1, A_2$ inside $\cQ_1$  for each $T$ there exists $t>T$
such that both sets $\varphi_t U_i \bigcap B_{R^*}(A_i) \bigcap
\cQ_j= \emptyset$ for all $i=1,2$ and $j=1,2$. In other words, the
image of $U_1$ (resp. $U_2$) at time $t$ bisects  the strip
$\cH_k$ into two parts and we can restrict this {\it nonempty}
intersection to not a long in the horizontal direction strip each.
Moreover, we shall need both intersections to be well separated in
the horizontal direction inside $\cH_k$, or dist$(A_1,A_2)=l \gg
R^*$.

Without the loss of generality one can assume that the points in
$\cH_k$ move in the direction from $A_1$ to $A_2$ under the vector
field $X_k$. Take two small constants $\epsilon, \tau>0$. Then
with positive probability the dynamics of $\varphi_{t,
t+\epsilon}$ is $C^0$-close to the dynamics of $\Phi^k_{1/\tau}$,
which is the time one map of the flow generated by $X_k/\tau.$
Define $\Lambda^k=B_{R^*}(A_1)\bigcap \cH_k$ and
$\cT^k_\tau=\Phi^k_{1/\tau} \Lambda^k.$ Then with positive
probability $\varphi_{t,t+\epsilon}(\cT^k_\tau)$ is $C^0$-close to
$\cT^k_\tau$. Moreover, $\varphi_{t+\epsilon}(U_1)$ divides it
into two parts. On the other hand, since $\varphi_{1/\tau}$ does
not move the set $\cQ_1$, we know that the set
$\varphi_{t+\epsilon}(U_2)$ still contains a point $q$ in
$\cQ_1\bigcap B_{R^*}(A_2).$ Now by Lemma \ref{CouldCr} there is a
positive probability that $W^s(q, t+\epsilon)$ intersects
$\cQ_2,$, moreover, since $R^*$ is sufficiently large,  the
segment of $W^s(q, t+\epsilon)$ between $q$ and the first crossing
of $\cQ_2$ lies inside $B_{R^*}(A_2)$. But if $l$ and $1/\tau$ are
much larger than $R^*$ this segment must cross the tongue
$\Phi_{1/\tau}(\Lambda).$ \qed

\subsection{Fences.} The following property of unstable foliations plays
the key role in our analysis. Let $\cS_N$ denote $N\times N$ square centered
at 0.

\begin{figure}[htbp]
  \begin{center}
    \begin{psfrags}
    \psfrag{H1}{$\cH_k^1$}
    \psfrag{H2}{$\cH_k^2$}
    \psfrag{W1}{$W^s(x_1, t)$}
    \psfrag{W2}{$W^s(x_2, t)$}
    \psfrag{Sn}{$S_N$}
    \psfrag{gt}{$\hat\gamma_t$}
            \includegraphics[height=3in, width= 4.5in,angle=0]{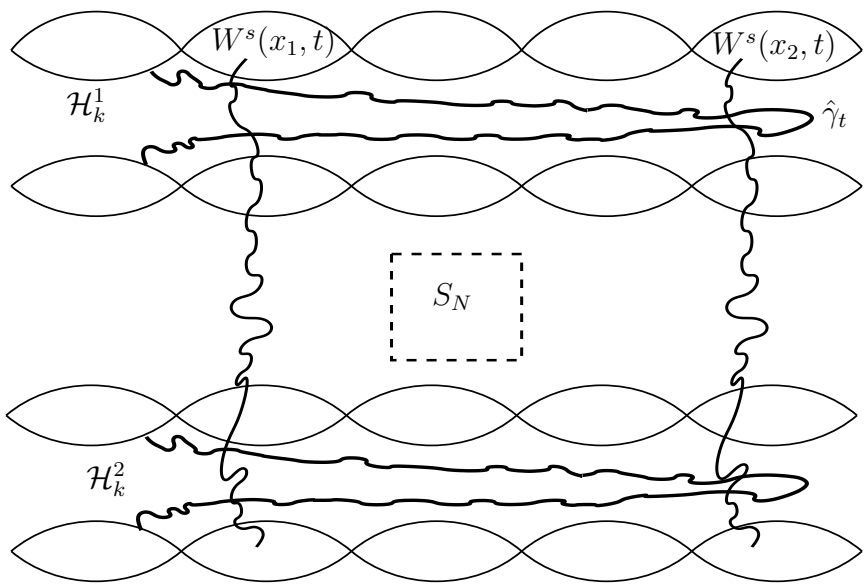}
    \end{psfrags}
    \caption{The construction used in the proof of transitivity of the stable
             lamination also gives the existence of fences.}
  \end{center}
\end{figure}

\begin{theorem} \label{contour}  For any curve
$ \gamma $  and any $N>0$ there almost surely exist a time $t$ and
contour $\Gamma$ consisting of two  pieces of  translated images
of $\gamma$ (${\gamma}_t+\vn_1$ and $\gamma_t+\vn_2),$ and two
stable curves, $W^s(x_1,t)$ and $W^s(x_2,t)$, such that
$\Int(\Gamma)$ contains  $\cS_N.$
\end{theorem}

\proof Let $X_k$ be the vector field with unbounded stream lines.
Let $\cH_k$ be a strip in $\R^2$, which contains some of the
unbounded stream lines of $X_k$. We can then find two parallel
translations of $\cH_k$ by integer vectors, to be denoted $\cH_k^1$
and $\cH_k^2$, which contain the square $\cS_N$  between them. Thus,
the plane is a union of five regions: the parallel strips $\cH_k^1$
and $\cH_k^2$, the strip between them, to be denoted by $\cR$, which
contains $\cS_N$, and two half planes, $\cP_1$, and $\cP_2$.

Let $\varepsilon > 0$ be fixed.  By Theorem \ref{WsTr}, for any $t
\geq 0$, there almost surely exists a stable manifold $W^s(x,t)$,
which intersects both half planes, $\cP_1$ and $\cP_2$. By
periodicity there are two stable manifolds, $W^s(x_1,t)$  and
$W^s(x_2,t)$, such that $\cS_N$ is contained in the part of $\cR$
between $W^s(x_1,t)$  and $W^s(x_2,t)$.

Let $r$ be such that ${\rm diam}(W^s(x_1,t) \bigcup W^s(x_2,t)
\bigcup \cS_N) \leq r$ with probability at least $1 -
\varepsilon$. Then by the hypoellipticity of the two point motion,
and the recurrence condition (Proposition \ref{LmP-PLD}), with
probability at least $1- 2 \varepsilon$ the curve $\gamma$
contains a part $ \widehat{\gamma}$, whose image  $
\widehat{\gamma}_t$ belongs to $\cH_k^1$, and intersects any curve
$\alpha$ connecting $\cP_1$ with $\cP_2$, with the property
${\rm diam} ( \alpha \bigcup \cS_N) \leq r$. We also require the same
of the integer translation of $ \widehat{\gamma}_t$ to $H_2$.

 The desired contour then consists
of $ \widehat{\gamma}_t$, its integer translation to $H_2$, and
the stable curves  $W^s(x_1,t)$ and $W^s(x_2,t)$. \qed

We say that $T$ is a {\it rearrangement} of $\cS_N$ if there is a
partition $T=\bigcup_j T_j$ of $T$ into finitely many pieces and a
bijection $\beta: T\to \cS_N$ such that on each $T_j$ we have
$\beta(x)=x+\vn_j,$ where $\vn_j\in\integers^2.$

From Theorem \ref{contour} it follows that there exists a contour
$\Gamma_0$  , which consists of two  pieces of  translations of
$\gamma$ (${\gamma}+\vn_1$ and $\gamma+\vn_2),$ and two stable
curves, $W^s(x_1,0)$ and $W^s(x_2,0)$, such that $\Int(\Gamma_0)$
contains a rearrangement of $\cS_N.$

\section{Upper Bound.}
\label{ub}
\subsection{Estimates in Probability.}
 \label{ssConvProb}
We first establish the asymptotic of the expectation of
$\tau^R(\gamma, tv)$.
\begin{lemma}
\label{ThProb} The following limit  $$ \lim_{t \rightarrow \infty}
\frac{\EXP\ \tau^R(\gamma,tv)}{t} = \|v\|^R $$ is uniform  in
$\gamma \in \cC.$
\end{lemma}

\proof Let $S=[-3R, 3R]^2.$ Then $\gamma\subset\Int(S)$ and so
$\tau^R(\gamma,tv)\geq \tau^R(S, tv).$ Hence it suffices to show
that $$ \frac{\EXP\ \tau^R(S,tv)}{t}\to \|v\|^R. $$ By Corollary
\ref{UI} it suffices to show that $$\left\{\tau^R(\gamma,
tv)-\tau^R(S,tv)\right\}_{\gamma\in\cC, t\geq 1, \|v\|=1}$$ is
tight. To prove the tightness let $G$ be the set of all
$\gamma$-fences containg a rearrangement of $S.$ For $\Gamma\in G$
let $S'(\Gamma)$ be some rearrangement of $S$ inside $\Gamma$
and let $\sigma(\Gamma)$ be the part of $\Gamma$ consisting of
stable manifolds.

Then for any $\varepsilon> 0$ there exists $r$ such that
\begin{equation}
\label{ShortFence}
\PROB\{\exists \Gamma\in G \text{ such that }
\diam(\Gamma \bigcup S)<r \;\;{\rm and}\;\;\sup_{t\geq 0 }
l(\sigma_t )<r\}\geq 1-\varepsilon~,
\end{equation}

Since by (\ref{ShortFence}) the distance between any two points,
first one from $S$, and the second one from $S'$, does not exceed
$r$, and since $S'$ is contained in $\Gamma$, we have
 that with probability not smaller than $1 - \varepsilon$
  for any $ A$ there is the inequality
$\tau^{r+R}(\Gamma,A)\leq \tau^R(S, A).$ Now
(\ref{ShortFence}) implies that $ {\rm dist}
(\gamma_{\tau^{r+R}(\Gamma,A)}, A)\leq 2r+R $ since $B(A,r+R)$
intersects either   $\gamma_{\tau^{r+R}(\Gamma,A)}$ or its integer
translation by distance no larger than $r$ , in which case
 $ {\rm dist}
(\gamma_{\tau^{r+R}(\Gamma,A)}, A)\leq r+R $, or $B(A,r+R)$
intersects the stable curve, and since the length of this curve is
at most $r$, there is a point in$\gamma_{\tau^{r+R}(\Gamma,A)}$ in
$B(A, 2r+R).$ Now by Proposition \ref{LmP-PLD} there is $t$ such that $$
\PROB\{\tau^R(\gamma_{\tau^{r+R}(\Gamma,A)}, A)>t\}\leq 2
\varepsilon.$$ This completes the proof. \qed

\begin{corollary}
\label{coro1}  For any $\gamma \in \cC$ we have $\lim_{t
\rightarrow \infty} \frac{\tau^R (\gamma,tv)}{t} = \|v\|^R $ in
probability.
\end{corollary}
\proof By Corollary \ref{UI} the family
$\{\frac{\tau^R(\gamma,tv)}{t}\}$ is tight. Let $\nu$ be some
limit distribution of this family. By Lemma \ref{Contains}
$\supp\nu\subset [\|v\|^R,\infty[.$ By Lemma \ref{ThProb} $\int t
d\nu(t)=\|v\|^R.$ Thus $\nu=\delta_{\|v\|^R}.$ \qed

\begin{corollary}
\label{AlmostContained} For any $\varepsilon>0$
$$
\lim_{t
\rightarrow \infty}  \PROB\{\cW^R_t\subset (1+\varepsilon) t\cB\} =
1~.
$$
\end{corollary}
\proof We need to show that
\beq \label{a1} \PROB\{\exists v:
\frac{\tau^R(\gamma,tv)}{t}<\frac{\|v\|^R}{1+\varepsilon}\}\to
0.
\eneq
Take $\varepsilon'\ll \varepsilon.$ Let $\{ v_j\}$ be an
$\varepsilon'$-net in $\partial\cB.$ Then by Corollary \ref{coro1}
\beq \label{b1}
\PROB\left\{\exists j:
\frac{\tau^R(\gamma,tv_j)}{t}<
\frac{\|v_j\|^R}{1+\frac{\varepsilon}{2}}\right\}\to0.
\eneq
However by Proposition \ref{LmP-PLD} for any $m$ for
sufficiently small $\varepsilon'$ we have
\begin{equation}
\label{c1}
\PROB\{\exists j:
\frac{\tau^R(\gamma,tv_j)}{t}<
\frac{\|v_j\|^R}{1+\frac{\varepsilon}{2}}|
\exists v:
\frac{\tau^R(\gamma,tv)}{t}<\frac{\|v\|^R}{1+\varepsilon}\}\geq 1 -
C_m t^{-m} .
\end{equation}
Now, (\ref{b1}) and (\ref{c1}) combined yield
(\ref{a1}). \qed
\subsection{Point-to-line passage.}
Given a line $l$ let $$\rho(l,t)=\sup_{\cC}
\frac{\EXP\ (\tau^R(\gamma, tl))}{t}.$$
The following results are proved similarly to a point-to-point case.
\begin{lemma} \label{PrP-L}
(a) $\frac{\EXP\ (\tau^R(\gamma, tl))}{t}$ converges
uniformly to $\rho(l)=\lim_{t\to\infty} \rho(l,t).$

(b) For any $ \varepsilon >0$ $$
\PROB\{\tau^R(\gamma,tl)>(\rho(l)+\varepsilon) t\}\leq C_N t^{-N}.$$
(c) $\frac{\tau^R(\gamma,tl)}{t}\to \rho(l)$ in probability.
\end{lemma}

The next lemma relates $\rho(l)$ to the norm $\|v\|^R$.

\begin{figure}[htbp]
  \begin{center}
    \begin{psfrags}
    \psfrag{S}{$S$}
    \psfrag{Pi}{$\Pi$}
    \psfrag{lp}{$l'$}
    \psfrag{tl}{$t l$}
            \includegraphics[height=4in, width= 4.5in,angle=0]{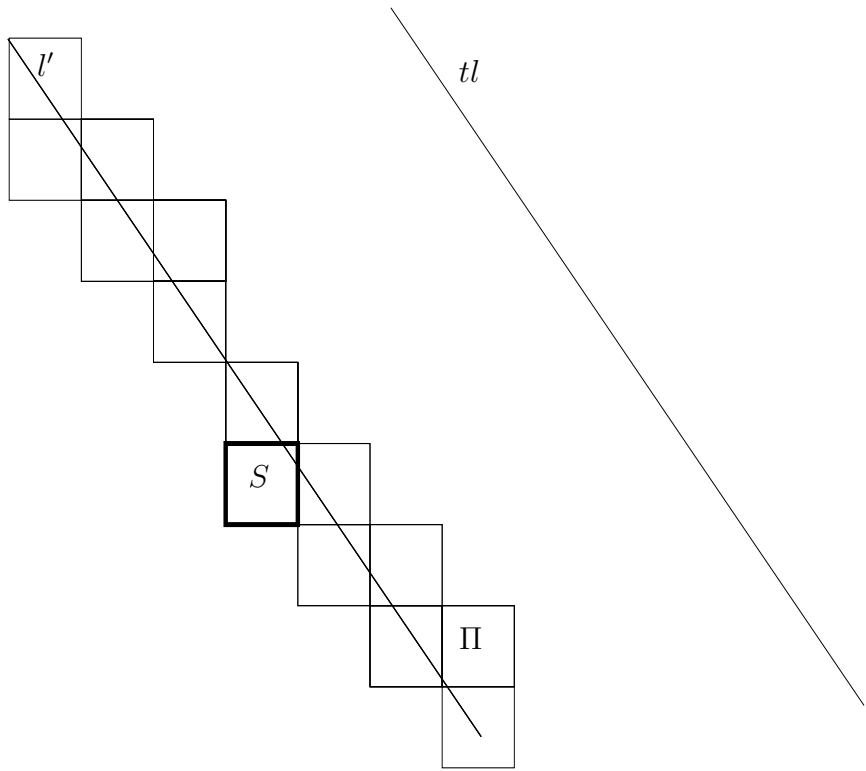}
    \end{psfrags}
    \caption{Line-to-line passage time differs little from square-to
line-passage time since all squares reach the lines at about the same time}
  \end{center}
\end{figure}

\begin{lemma} For any line $l$
$$ \rho(l)=\inf_{v\in l} \|v\|^R . $$
\end{lemma}
\proof This follows from Lemma \ref{PrP-L}, Lemma \ref{Contains}
and Corollary \ref{AlmostContained}. \qed

 Consider now
$$\rho^*(l)=\lim_{t\to\infty} \frac{\EXP\ (\tau^R(l', tl))}{t}$$
where $l'$ is the line parallel to $l$ passing through origin.
\begin{lemma} For any line $l$ we have the equality
$$ \rho^*(l)=\rho(l).$$
\end{lemma}
\proof Let $\Pi$ be the union of fundamental domains containing
$l'.$ Applying the reasoning of subsection \ref{ssConvProb} we get
$$
\EXP\ \left(\frac{\tau^R(l',tl)-\tau^R(\Pi, tl)}{t}\right)\to
0.
$$
Let $S$ be one fundamental domain from $\Pi.$ Then$$
\EXP\ \left(\frac{\tau^R(S,tl)-\tau^R(\Pi, tl)}{t}\right)\to 0.$$
Subtracting we get $$ \EXP\ \left(\frac{\tau^R(l',tl)-\tau^R(S,
tl)}{t}\right)\to 0.$$ Passing to the limit as $t\to\infty$ we get
the statement required. \qed

\subsection{Almost sure convergence.}
\begin{theorem} For any $\varepsilon > 0$
almost surely $\cW^R_t\subset (1+\varepsilon) t \cB$ for large
enough $t.$
\end{theorem}
\proof Again choose some $\varepsilon_1\ll\varepsilon$ and let
$\{v_j\}$ be an $\varepsilon_1$-net on $\partial\cB.$ Let
$\cB_{\varepsilon_1}$ be the region bounded by support lines of
$\cB$ passing through $\{v_j\}.$ It suffices to prove that almost
surely $\cW^R_t\subset t(1+\varepsilon) \cB_{\varepsilon_1} $ for
large $t.$ But this accounts to showing that almost surely
\begin{equation}
\label{LnNotTooFast} \tau^R(\gamma, tl)\geq (1-\varepsilon) t
\rho(l)
\end{equation}
for large $t.$
Let $t^*$ be such that
\begin{equation}
\label{ReachLine}
\EXP\ (\tau^R(l',tl))>(1-\frac{\varepsilon}{2})(t+1)\rho(l)
\end{equation}
for $t>t^*.$ Let $l_j=j t^* l.$ Let $\Delta_j$ be the first moment
after $\tau^R(\gamma, l_j)$ such that
$l_j(\tau^R(\gamma,l_j)+\Delta)\bigcap l_{j+1}\neq\emptyset.$ It
follows that $\tau^R(\gamma, l_j)\geq \sum_{k=0}^{j-1} \Delta_k.$
Thus (\ref{LnNotTooFast}) follows by (\ref{ReachLine}) and Lemma
\ref{SMartLD}. \qed

\end{document}